\documentclass[12pt]{article}
\usepackage{verbatim}
\usepackage{amssymb}
\usepackage{amsthm}
\usepackage{amsmath}
\usepackage{graphicx}
\usepackage{mathrsfs}
\AtBeginDocument{%
 \def\MR#1{} 
}
\usepackage{tikz}
\usepackage{tikz-cd}  
\usepackage{yfonts}
\usepackage{makeidx}
\usepackage{enumerate}
\numberwithin{equation}{section}

\usepackage{color}
\usepackage[backrefs]{amsrefs}
\usepackage[pagebackref,hypertexnames=true, colorlinks, citecolor=black, linkcolor=blue, urlcolor=red]{hyperref}

\newcommand\CC{{\mathcal C}}
\renewcommand\SS{{\mathscr S}}
\newcommand\h{\CC_\h}
\renewcommand\O{\Omega}
\newcommand\cbhd{CB(\h)^d}

\newcommand\hinf{H^\infty}
\def\norm#1{\| #1 \|}
\newcommand\BB{{\mathcal B}}

\usepackage{latexsym}
\usepackage{amssymb}
\usepackage{euscript}

\let\cal=\mathcal      

\def\mcc{M\raise.5ex\hbox{c}C}
\def\mccarthy{M\raise.5ex\hbox{c}Carthy}

\def\h{{\cal H}}

\def\K{{\cal K}}
\def\M{{\cal M}}

\def\l{\lambda}

\def\la{\langle}
\def\ra{\rangle}
\def\={\ = \ }

\def\A{{\cal A}}
\def\BB{{\cal B}}
\def\CC{{\cal C}}

\def\PP{{\cal P}}

\def\DD{D_\delta(k,w)}

\def\C{\mathbb C}

\def\D{\mathbb D}

\def\be{\setcounter{equation}{\value{theorem}} \begin{equation}}
\def\ee{\end{equation} \addtocounter{theorem}{1}}
\def\beq{\begin{eqnarray*}}
\def\eeq{\end{eqnarray*}}
\def\se{\setcounter{equation}{\value{theorem}}} 
\def\att{\addtocounter{theorem}{1}}
\def\vs{\vskip 5pt}

\def\bp{{\sc Proof: }}
\def\ep{{}{\hfill $\Box$} \vskip 5pt \par}

\def\bl{\begin{lemma}}
\def\el{\end{lemma}}
\def\bt{\begin{theorem}}
\def\et{\end{theorem}}
\def\bprop{\begin{prop}}
\def\eprop{\end{prop}}
\def\bd{\begin{definition}}
\def\ed{\end{definition}}
\def\br{\begin{remark}}
\def\er{\end{remark}}
\def\bexer{\begin{exercise}}
\def\eexer{\end{exercise}}

\newtheorem{theorem}{Theorem}[section]
\newtheorem{prop}[theorem]{Proposition}
\newtheorem{lemma}[theorem]{Lemma}

\newtheorem{definition}[theorem]{Definition}

\numberwithin{equation}{section}

\renewcommand\DD{{\mathcal D}}
\newcommand\hig{H^\infty_{\rm gen}}

\title{
 Calcular algebras}
\author{Jim Agler
\thanks{Partially supported by National Science Foundation Grant
 DMS 1665260}
\and
John E. M\raise.5ex\hbox{c}Carthy
\thanks{Partially supported by National Science Foundation Grant  
DMS 1565243
}
\and
N. J. Young\thanks{Partially supported
by UK Engineering and Physical Sciences Research Council grants EP/K50340X/1 and  EP/N03242X/1, and 
London Mathematical Society grants 41219 and 41730,
{\bf  MSC} [2010]:  15A54, 32A99, 58A05, 58J42
}}

\begin{document}
\maketitle

\centerline{To the memory of Richard Timoney}

\begin{abstract}
A calcular algebra is a subalgebra of $H^\infty(\Omega)$ with norm given by
$\| \phi \| = \sup \| \phi(T) \|$ as $T$ ranges over a given class of commutative $d$-tuples of operators with Taylor spectrum in $\O$. We discuss what algebras arise this way, and how they can be represented.
\end{abstract}
\section{Introduction}

Let $\Omega$ be a bounded open set in $\C^d$. 
We say that a class $\CC$ is {\em subordinate to $\O$} if:
\begin{enumerate}
\item[(i)]
Each element $T$ of $\CC$
is a commuting $d$-tuple of bounded operators on a Hilbert space, with 
its Taylor spectrum\footnote{For a definition of Taylor spectrum of a commuting tuple, see \cite{tay70a}.}
 $\sigma(T)$  in $\Omega$.
\item[(ii)]
For some non-zero Hilbert space $\h$, 
the set of scalars
\[
\{ (\l^1, \dots, \l^d ) : \l \in \O \} \ \subseteq \ \CC ,
\]
where we think of $\l$ as a $d$-tuple of scalar multiples of the identity  acting on $\h$.
(Note that we use superscripts to denote the coordinates.)
\end{enumerate}

Given a class $\CC$ subordinate to $\Omega$, we define 
$H^\infty(\CC)$ to be those holomorphic functions on
 $\O$ for which
 \[
 \| \phi \|_{\CC} \= \sup \{ \| \phi(T) \| : T \in \CC \}
 \]
 is finite. It can be shown (see Prop. \ref{propb0} below) that this algebra is always complete, so it  is a Banach algebra, which by Property (ii) is always contained contractively in
 the algebra $H^\infty(\O)$ of bounded holomorphic functions on $\O$.
(We are using $H^\infty$ in two apparently different ways, but identifying $\Omega$ with the
set of scalars makes the two usages agree).
Any Banach algebra of holomorphic functions arising in this way we shall call a
{\em calcular algebra over $\Omega$}.

 We shall call the closed unit ball 
 of $H^\infty(\CC)$ the  {\em Schur class of $\CC$}, and denote it by 
  $\SS(\CC)$.
 \[
 \SS(\CC) \= \{ \phi \in {\rm Hol}(\O) : \| \phi(T) \| \leq 1, \ \forall \ T \in \CC \}.
 \]
Let $\SS(\Omega)$ denote the closed unit ball
 of  $H^\infty(\O)$.
 
If $\h$ is  a Hilbert space (we shall always assume that Hilbert spaces are  not zero-dimensional to avoid
trivialities), 
let $\cbhd$ denote the set of commuting $d$-tuples of elements of $B(\h)$, the bounded linear operators on $\h$.
Given a set $S$ of bounded holomorphic functions on $\O$, and a Hilbert space $\h$, one can form the set
 \[
 \h(S) \= \{ T  \in \cbhd :\ \sigma(T) \subseteq \O \ {\rm and\ } \| \phi(T) \| \leq 1 \ \forall \ \phi \in S \}.
 \]
 
 If $\h$ is a Hilbert space, $\CC \subseteq \cbhd$ and $S \subseteq \SS(\Omega)$, then tautologically one has
 \se\att
 \begin{eqnarray}
 \label{eqa09}
 \h(\SS(\CC)) & \ \supseteq  \ & \CC\qquad
 \\
  \label{eqa10}
 \SS(\h(S)) & \supseteq & S   .
 \end{eqnarray}
\att
 Typically these inclusions will be strict. For example, let $d=1$, and let $\Omega$ be the open unit disk $\D$.
 Let $\h$ be any Hilbert space, and let $\CC$ be the set $\{ \lambda I : \lambda \in \D \}$.
 Then $\SS(\CC)$ will equal $\SS(\D)$, and, by von Neumann's inequality \cite{vonN51},
 $\h(\SS(\D))$ will consist of all contractions on $\h$ whose spectrum is in $\D$.
 Likewise if $S$ just contains the function $z$, then $\h(S)$ will be the contractions on $\h$ whose spectrum is in $\D$,
 and the Schur class of this set will be all of $\SS(\D)$.
 
 Our first result is that the operations $\h$ and $\SS$ stabilize after 3 steps, provided $\h$ is infinite dimensional.
 
 Notation: If $T$ is a commuting $d$-tuple of bounded operators on a Hilbert space $\h$, we call $\h$ the carrier of
 $T$, and write $ \h = {\rm car}(T)$.
 
 \bt
 \label{thm1}
 Let $\O$ be a bounded open set in $\C^d$, and let $\CC$ be any class subordinate
 to $\O$. Let $S$ be a non-empty subset of $\SS(\Omega)$.
For any Hilbert space $\h$, we have
\be
 \label{eqa20}
  \h(\SS(\h(S))) = \h(S) .
  \ee
 If the dimension of $\h$ is either infinite, or greater than or equal to 
 $\sup \{ {\rm dim} ({\rm car} (T)) : T \in \CC \}$,
 then
 \be
 \SS(\h(\SS(\CC))) \= \SS(\CC) 
 \label{eqa30}
 \ee
 \et
 \bp
By \eqref{eqa10}, we have
\be
 \label{eqa40}
 \h(\SS(\h(S))) \subseteq  \h(S) .
 \ee
Suppose now that $T \in \h(S)$, and $\phi$ is any function in $\SS(\h(S))$.
Then $\| \phi(T) \| \leq 1$, so $T$ is in $ \h(\SS(\h(S)))$, proving \eqref{eqa20}.

By \eqref{eqa10} again, with $S = \SS(\CC)$, we get
\be
\label{eqa50}
\SS(\CC) \ \subseteq \ 
 \SS(\h(\SS(\CC))) .
 \ee
Now, assume the dimension of $\h$ is as in the hypothesis.
Let $\phi \in \SS(\h(\SS(\CC)))$, and let $T \in \CC$, with ${\rm car}(T) = \K$.
We need to show $\| \phi(T) \| \leq 1$.
 If the dimension of $\h$ is  equal to the dimension of $\K$, then $T$ is unitarily equivalent to a $d$-tuple
 $R$ on $\h$, and $R \in \h(\SS(\CC))$ since $\| \psi(R) \| = \| \psi(T) \| \leq 1$ for every $\psi $ in
 $\SS(\CC)$. Therefore $\| \phi (T) \| = \| \phi (R) \| \leq 1$, and we are done.

If the dimension of $\h$ is larger than the dimension of $\K$, write $\h = \h_1 \oplus \h_2$
 where ${\rm dim}(\h_1) = {\rm dim}(\K)$, and let $R_1$ on $\h_1$ be unitarily equivalent to $T$.
 Choose $\lambda = (\lambda^1, \dots, \lambda^d)$ in $\Omega$, and let 
 \[
 R \= (R_1^1 \oplus \lambda^1 I_{\h_2} , \dots, R_1^d \oplus \lambda^d {I}_{\h_2} ) ,
 \]
 Then for any $\psi \in O(\Omega)$, the set of holmorphic functions on $\Omega$,
  we have $\psi(R) = \psi(R_1) \oplus \psi(\lambda) I_{\h_2}$, 
 so if $\psi \in \SS(\C)$, we have $\| \psi(R) \| \leq 1$, and therefore $R \in \h(\SS(\CC))$. Now we get
 $\| \phi(T) \| \leq \| \phi(R) \| \leq 1$, and again we are done.
 
 Finally we consider the case where $\h$ is infinite dimensional, but the carriers of the elements of $\CC$ 
 may have larger dimension. We can assume without loss of generality that $\h$ is separable. 
We need to find $R \in  \h(\SS(\CC))$ with $\| \phi (R) \| = \| \phi(T) \|$.
To do this, it is sufficient to show that there is a separable subspace $\K_1$ of $\K$ that is reducing
for $f(T)$ for every $f \in O(\O)$ and such that $\| \phi (T ) |_{\K_1} \| = \| \phi(T) \|$; for then we can choose
$R_1$ on $\h$ unitarily equivalent to $ P_{\K_1} T |_{\K_1} $, where $P_{\K_1}$ is projection onto $\K_1$;
the fact that $\K_1$ is reducing means that $ \phi(  P_{\K_1} T |_{\K_1} ) =  P_{\K_1} \phi( T) |_{\K_1}$.

 Observe that $\{ f(T) : f \in O(\O) \}$ has a countable dense subset $\DD$ in the norm topology of $CB(\K)^d$, since $O(\O)$ is separable.
Let $u_j$ be  a sequence of unit vectors in $\K$  such that $\| \phi(T) u_j \| \to \| \phi(T) \|$.
Let $\K_1$ be the closed linear span of finite products of elements of $\DD \cup \DD^*$ applied to
finite linear combinations of the vectors $u_j$.
By $\DD^*$ we mean
\[
\DD^* \ = \ \{ ((T^1)^*, \dots , (T^d)^*) \ : \ (T^1, \dots, T^d) \in \DD \}.
\]
 Then $\K_1$  is a separable subspace of $\K$ on which $\phi(T)$
achieves its norm and that is reducing for every $f(T)$. \ep

For a given class $\CC$, it is of interest to know the smallest dimension of $\h$ that gives equality
in 
\eqref{eqa30}.  

\begin{exam}
Let $\Omega = \D$, and let $\CC$ be all contractions with spectrum in $\D$.
Then we can take $\h$ to be one dimensional. 
Similarly, if $\Omega = \D^2$, and $\CC$ is all pairs of commuting contractions with
spectrum in $\D^2$, And\^o's inequality \cite{and63} yields that we can take $\h$ to be one dimensional again.

However, if $d \geq 3$, we let $\Omega = \D^d$, and $\CC$ be the class of all $d$-tuples of commuting contractions
with spectrum contained in $\D^d$, then $\SS(\CC)$ is the Schur-Agler class, a proper subset of $\SS(\D^d)$
\cite{var74, cradav}. If $\h = \C^n$, then 
$\h (\SS(\CC))$ will be all $d$-tuples of commuting contractive $n$-by-$n$ matrices with spectrum in $\D^d$.
In \cite{kn16}, it is shown that if $n=3$, then
\[
\SS(\C^3(\SS(\CC))) \= \SS(\D^d) .
\]
It is unknown what the minimal dimension of $\h$ must be in this case to get equality in
\eqref{eqa30}, or even whether it must be infinite.
\end{exam}
\begin{exam}
Let $\K$ be a Hilbert function space on $\Omega$ with reproducing kernel $k$.
The multiplier algebra ${\rm Mult}(\K)$ is always a calcular algebra.
Indeed, for each finite set $F = \{\lambda_1, \dots, \l_n \} \subset \O$, let $T_F$
be the commuting $d$-tuple $(T_F^1, \dots , T_F^d)$ acting on the $n$-dimensional subspace
of $\K$ spanned by the kernel functions $\{ k_{\lambda_j} : 1 \leq j \leq n \}$ defined by 
\[
T_F^r k_{\lambda_j} \ = \ \overline{\lambda_j^r} k_{\l_j} \qquad 
1\leq r \leq d, 1 \leq j \leq n.
\]
Define 
\[
\CC = \{ T_F^* : F \ {\rm a\ finite\ subset\ of \ } \O \}.
\]
It is straightforward to show that $H^\infty(\CC) = {\rm Mult}(\K)$.
\end{exam}
Many other examples of calcular algebras are given in 
 \cite[Chapter 9]{OpAn}.

\section{When is a Banach algebra a calcular algebra?}
\label{secb}

\begin{prop}
\label{propb0}
Let $\CC$ be subordinate to $\O$. Then $H^\infty(\CC)$ is a Banach algebra.
\end{prop}
\bp
We need to prove completeness.
Consider a Cauchy sequence $\{\phi_n\}$ in $\hinf(\CC)$.
Since $\CC$ is subordinate to $\O$,
 $\{\phi_n\}$ is a Cauchy sequence in $\hinf(\Omega)$. Therefore, as $\hinf(\Omega)$ is complete, there exists $\phi \in \hinf(\Omega)$ such that
\be\label{nor40}
\sup_{\lambda \in \Omega}\ |\phi_n(\lambda)-\phi(\lambda)| \to 0\  \text{ as }\  n \to \infty.
\ee
We claim that
\be\label{nor50}
\phi \in \hinf(\CC)
\ee
and
\be\label{nor60}
\phi_n \to \phi\  \text{ in }\ \hinf(\CC).
\ee

To prove statement \eqref{nor50}, note that for each $T\in \CC$,  we have $\Omega$ is a neighborhood of $\sigma(T)$. 
Consequently, 
 continuity of the functional calculus implies that
\be\label{nor70}
 \phi_n(T) \to \phi(T) \qquad \ \forall\  T\in \CC.
\ee
Also, as $\{\phi_n\}$ is a Cauchy sequence in $\hinf(\CC)$, there exists a constant $M$ such that
\[
 \norm{\phi_n}_\CC \le M \qquad \forall\  n.
\]
Therefore, if $T\in \CC$,
\[
\norm{\phi(T)}_\CC\;=\;\lim_{n \to \infty}\norm{\phi_n (T)}\le \limsup_{n \to \infty} \norm{\phi_n}_\CC \;\le \; M.
\]
Hence,
\[
\norm{\phi}_\CC\; = \;\sup_{T\in\CC}\norm{\phi(T)}_\CC\; \le \; M,
\]
which proves the membership \eqref{nor50}.

To prove the limiting relation \eqref{nor60}, let $\varepsilon>0$. Choose $N$ such that
\[
m,n\ge N \implies \norm{\phi_n-\phi_m}_\CC<\varepsilon.
\]
By definition of the norm, this means
\[
 m,n\ge N \implies \norm{\phi_n(T)-\phi_m(T)}<\varepsilon \qquad \forall\ T \in \CC.
\]
Letting $m \to \infty$ and using statement  \eqref{nor70} we deduce that
\[
n\ge N \implies \norm{\phi_n(T)-\phi(T)}<\varepsilon\qquad \forall\ T\in \CC.
\]
Hence, since $\norm{\phi_n-\phi}_\CC = \sup_{T \in\CC}\norm{\phi_n(T)-\phi(T)}$,
\[
n\ge N \implies \norm{\phi_n-\phi}_\CC\le\varepsilon.
\]
\ep

Let $\A$ be a unital Banach algebra contractively contained in $H^\infty(\Omega)$.
When can it be realized as a calcular algebra? Let $S$ be its unit ball. By Theorem \ref{thm1},
$\A$ is a calcular algebra if and only if $\SS(\h(S)) = S$, where $\h$ is an infinite dimensional Hilbert space.

This imposes a constraint on $\A$. In particular, 
 there must be an isometric homomorphsim from $\A$ into $B(\K)$ for some Hilbert space $\K$.
 There is another constraint which stems from the requirement that all the operators in the class
 have spectrum in the open set $\Omega$.
\begin{prop}
\label{propb1}
If $\A$ is a calcular algebra, then:

{\rm (i)}  There is an isometric homomorphism into $B(\K)$
for some Hilbert space $\K$. 

{\rm  (ii)} If $\phi_n$ is a bounded sequence in $\A$ that converges uniformly on compact subsets
of $\O$ to a function $\psi$, then $\psi \in \A$, and $\| \psi \| \leq \liminf \| \phi_n \|$.

\end{prop}
\bp (i)
Suppose $\A$ is $H^\infty(\CC)$ for some class $\CC$ subordinate to an open set $\Omega$.
Let $\h$ be any infinite dimensional Hilbert space, and let $\CC_1 = \h(\SS(\CC))$.
By Theorem~\ref{thm1}, we have
\be
\label{eqb1}
\A \= H^\infty(\CC_1) .
\ee
Let $\K$ be
the direct sum of cardinality($\CC_1$) copies of $\h$, with the sum indexed by $\CC_1$.
Define a map $\pi : \A \to B(\K)$ by
\[
\pi(\phi) \= \oplus_{T \in \CC_1} \phi(T) .
\]
Then $\pi$ is a homomorphism, and by \eqref{eqb1} it is isometric.

(ii) Let $\phi_n$ be a bounded sequence in $H^\infty(\CC)$ converging to $\psi$ locally uniformly on $\O$.
Without loss of generality, we may assume that each $\phi_n$ is in $\SS(\CC)$.
For each $T$ in $\CC$, since $\sigma(T) \subseteq \Omega$, it follows from the continuity of
the functional calculus that $\psi(T)$ is the limit in norm of $ \phi_n(T)$, so $\psi$
is in $\SS(\CC)$. 
Replacing $\phi_n$ by a subsequence whose norms converge to $ \liminf \| \phi_n \|$ gives the last inequality.
\ep
Remark: If one defines $ S =  \oplus_{T \in \CC_1} (T)$, then one can interpret $\phi(S)$ 
as $\pi(\phi)$. However, the spectrum of $S$ will be $\overline{\Omega}$, so $ S$
is not contained in any class subordinate to $\Omega$. The Taylor functional calculus is
 defined only for functions holomorphic on a neighborhood of the Taylor spectrum
of the $d$-tuple.

\vs

A necessary condition for a Banach algebra to be isometrically isomorphic to an algebra of operators
on a Hilbert space is that it satisfies
von Neumann's inequality: $\| p(x) \| \leq \| p \|_{H^\infty(\D)}$ for any $x$ in the unit ball of
the Banach algebra, and any polynomial $p$. It is not known whether this condition is sufficient.

Calcular algebras come with a sequence of matrix norms. If $[ \phi_{ij} ]$ is an $n$-by-$n$ matrix
of elements of $H^\infty(\CC)$, one can define 
\[
\| [ \phi_{ij} ] \|_n \= \sup \{ \| [\phi_{ij}(T)] \| : T \in \CC \} ,
\]
where the norm on the right-hand side is the operator norm on ${\rm car}(T) \otimes \C^n $.
By a similar argument to Proposition~\ref{propb1}, one can show that calcular algebras
have completely isometric homomorphic embeddings into some $B(\K)$.

Algebras that can be completely isometrically realized in this way are characterized
by the Blecher-Ruan-Sinclair theorem \cite{brs90}, \cite[Cor. 16.7]{pau02}.
This says that the algebra $\A$ must satisfy the Ruan axioms:
\beq
\forall n \in {\mathbb N}, \forall  a \in M_n(\A), \forall  X,Y \in M_n(\C), \quad 
\| XaY \|_n &\ \leq \  &\| X \| \| a \|_n \| Y \| \\
\forall m,n \in {\mathbb N}, \forall  a \in M_m(\A), b \in M_n(\A), \quad
\| a \oplus b \|_{m+n} &=& \max{ \| a \|_m, \| b \|_n} ,
\eeq
 and hence be isometrically realizable as an
operator space; and the matrix multiplication at each level $n$ must be contractive,
{\em i.e.} if $a = [a_{ij} ]$ and $b = [b_{ij} ]$ are in $M_n(\A)$,
then
\[
\| [ \sum_{k=1}^n a_{ik} b_{kj} ] \|_n 
\ \leq \
\| [ a_{ij}] \|_n \| [ b_{ij} ] \|_n  .
\]
It is straightforward to check that a calcular algebra satisfies the hypotheses of the Blecher-Ruan-Sinclair theorem.

We do not know in general what intrinsic
 necessary and sufficient conditions on a sub-algebra of $\hinf(\O)$
make it a calcular algebra; we can say something 
with strong convexity assumptions.
 Let $\PP = \C[z_1, \dots, z_d]$ denote the polynomials.
If $f$ is a function and $r > 0$, define $f_r$ by
$f_r(z) = f(rz)$.

\bt
\label{thmb2}
Let $\A$ be a unital Banach algebra that is contractively contained in $\hinf(\O)$, for some bounded open convex set $\O$ in $\C^d$ that contains $0$.
Suppose that $\PP$ is contained in $\A$ and that for every function $\phi \in \A$, there is a
sequence  in $\PP$ that is bounded in norm by $\| \phi \|$
 and converges to $\phi$ locally uniformly
on $\O$. Suppose moreover that for every polynomial $p \in \PP$,
we have  $\| p_r \| \leq \| p \|$ for $0 < r < 1$.

Then $\A$ is a calcular algebra over $\Omega$ if and only if the necessary conditions
of Proposition \ref{propb1} hold.
\et
\bp
Suppose both conditions hold, and $\pi$ embeds $\A$ isometrically in $B(\K)$.
For each of the coordinate functions $z^j, 1 \leq j \leq d,$ define
$T^j = \pi(z^j)$. Let $T \in CB(\K)^d$ be the tuple $(T^1, \dots, T^d)$.
Then for any polynomial $p \in \PP$ we have
$\pi(p) = p(T)$.
Moreover, if $p$ has no zeroes on $\overline{\O}$, then
\[
\pi( p \  \frac{1}{p} ) \= 1_\K \= p(T) \pi( \frac{1}{p} ) ,\]
so
\[
 \pi( \frac{1}{p} ) \= p(T)^{-1} .
 \]
As $p$ ranges over affine functions whose zero sets are hyperplanes not intersecting $\overline{\O}$,
we see that $\sigma(T)$ must be contained in $\overline{\O}$.

We want the elements of $\CC$
to have spectrum in $\O$.
 Let $\CC = \{  r T  : 0 \leq r < 1 \}$.

For any polynomial $p$ and any sequence $r_n \uparrow 1$ we have
\beq
\| p \|_\A & \ = \ & 
\| \pi(p ) \| \\
&=& \| p(T) \| \\
&=& \lim_{n \to \infty} \| p (r_n T) \|  \\
&\leq & \| p \|_\CC  \\
&=& \sup_{0 < r < 1} \| p_r (T) \| \\
&=& \sup_{0 < r < 1} \| \pi(p_r ) \| \\
&\leq & \| p \|_\A  .
\eeq
So $\A$ and $H^\infty(C)$ assign the same norm to polynomials.

Let $\psi$ be in  $\A$ of norm $1$. By hypothesis, there is a  sequence of polynomials
$q_n$ that converges locally uniformly to $\psi$, with $\| q_n \|_\A \leq 1$. Therefore 
for each $ 0 \leq r < 1$, \[
\| \psi(rT) \| \= \lim_{ n \to \infty}\|  q_n(rT) \| \ \leq \ 1 .
\]
Therefore $\psi$ is in the unit ball of $H^\infty(C)$, and hence $\A$  is
contractively contained in $ H^\infty(C)$.

Conversely, let $\phi \in \SS(C)$. Since $\O$ is convex, $\phi_r$ will converge
to $\phi$ locally uniformly on $\O$ as $r \uparrow 1$.
Fix $r < 1$. There is a sequence $q_n$ of polynomials that converges uniformly to $\phi_r$
on a neighborhood of $\overline{\O}$. Therefore $\lim_{n \to \infty} q_n(T) = \phi_r(T)$
is a contraction, and so by Property (ii) we have
\[
\| \phi_r \|_\A \ \leq \ 1 .
\]
By a diagonalization argument, we can modify this construction to find polynomials $q_n$
in the unit ball of $\A$ that converge locally uniformly to $\phi$, and hence
\[
\| \phi \|_\A \leq 1 .
\]
So  $H^\infty(C)$ is contractively contained in $\A$,
and hence the two algebras are isometrically isomorphic.
\ep

\begin{exam}
The disk algebra $A(\D)$ cannot be a calcular algebra, since it fails (ii).
However, there are subalgebras of the disk algebra that are the multiplier algebra 
of some Hilbert function spaces on the disk, {\em e.g.} the space with reproducing kernel
\[
k(w,z ) \= \sum_{n=0}^\infty (n+1)^2 z^n \bar w^n .
\]
Multiplier algebras for spaces of holomorphic functions are always calcular,
as shown in Example 1.9.
\end{exam}

\begin{prob} Find necessary and sufficient conditions for a subalgebra of $H^\infty(\O)$ to
be a calcular algebra.
\end{prob}
\section{Realization formulas}

In \cite{dmm07} and \cite{dm07}, Dritschel,  Marcantognini,  and  McCullough
proved a very general realization formula, building on 
work of Ambrozie and Timotin in 
\cite{at03},
which can be adapted to our current setting.

Let $S$ be a set of  functions from a set $X$ to the   unit disk ${\D}$.
In this section, we shall make the standing assumption that $S$ restricted to any finite set $F$ generates, as an algebra,
all the complex-valued functions on $F$. 

We define $K_S$ to be the set of kernels on $X$ that satisfy
\[
K_S \= \{ k \ |\  ( 1 - \psi(z) \bar \psi(w) ) k(z,w) \ \geq \ 0 \quad \forall\ \psi \in S \} .
\]
We define $A^\infty(K_S)$ to be
\[
A^\infty(K_S) = \{ \phi: X \to \C \ | \ \exists M  \geq 0 \ {\rm s.t. }\
(M^2 -  \phi(z) \bar \phi(w) ) k(z,w) \ \geq \ 0 \ \forall\ k \in K_S \} ,
\]
and define $\| \phi \|$ to be the smallest $M$ that works.

Endow $S$ with the topology of pointwise convergence.
 Let $C_b(S)$
denote the continuous bounded functions on $S$, which we think of as a C*-algebra. Let $E : X \to C_b(S)$ be the 
evaluation map $E(z) (\psi) = \psi(z)$, and let $E(w)^*$ mean the complex conjugate of this, the adjoint in the C*-algebra, 
$E(w)^* (\psi) = \overline{\psi(w)}$.

If $\psi$ is a function from $X$ to $\C$, we say it has a {\em network realization formula}
if there exists a Hilbert space $\M$, a unital *-representation $\rho : C_b(S) \to B(\M)$, and a unitary
$U : \C \oplus \M \to \C \oplus \M$ that in block matrix form is
\[
U \= \begin{pmatrix} A&B \\C&D \end{pmatrix}
\]
so that
\be
\label{eqc1}
\psi(z) \= A + B \rho(E(z)) ( I - D \rho(E(z)))^{-1} C.
\ee

If $\BB$ is a C*-algebra, a positive kernel on a set $X$ with values in $\BB^*$, the dual of $\BB$,
is a function $\Gamma : X \times X \to \BB^*$ such that for every finite
 set $F \subset X$, and every $f : F \to \BB$ we have
\[
\sum_{z,w \in F} \Gamma(z,w)(f(w)^* f(z) ) \geq 0 .
\]
Here is the Dritschel,  Marcantognini,  and  McCullough theorem.
\bt
\label{thmc1}
Let $S$ be a set of  functions from $X$ to  ${\D}$,
and let $\phi : X \to \overline{\D}$. The following are equivalent:
\begin{itemize}
\item[\rm (i)]
$\phi \in A^\infty(K_S)$ and $ \| \phi \|_{A^\infty(K_S)} \leq 1 $.
\item
[\rm (ii)]
 For each finite set $F \subseteq X$ there
exists a positive kernel $\Gamma : F \times F \to C_b(S)^*$ so that,
for all $z,w \in F$,
\be
\label{eqc7}
1 - \phi(z) \overline{\phi(w)} \= \Gamma(z,w) ( 1 - E(z) E(w)^*) .
\ee
\item
[\rm (iii)] $\phi$ has a network realization formula.
\end{itemize}
\et

Now let us assume that the functions in $S$ are all holomorphic functions on the open set $\O$ in $\C^d$.
By definition, we always have $S$ is contained in the unit ball of $A^\infty(K_S)$, so
when $\h$ is infinite dimensional we  have $H^\infty(\h(S))$ is contractively contained
in  $A^\infty(K_S)$ by Theorem~\ref{thm1}.
We shall show in Theorem~\ref{thmc2} and Proposition \ref{propc3} that the converse holds if $S$ is finite, or if a certain generic assumption holds.

We shall say that $T$ is a {\em generic matrix $d$-tuple on $\Omega$} if,
for some $n \in {\mathbb N}$, we have that  $T$ is a $d$-tuple of commuting
$n$-by-$n$ matrices that have a common set of $n$ linearly independent eigenvectors with
distinct joint eigenvalues, which means there are $n$  linearly independent eigenvectors $v_j$
in $\C^n$ so that 
\be
\label{eqc3}
T^r v_j \= \lambda_j^r v_j, \qquad 1\leq r \leq d,\ 1 \leq j \leq n ,
\ee
and the $n$ points $\lambda_j = (\lambda_j^1, \dots, \lambda_j^d)$ are distinct
points in $\Omega$. The advantages of working with generic $d$-tuples were pointed out in
\cite{amhirosh}.

We shall define an algebra $\hig(S)$ to be the holomorphic functions on $\O$
for which the norm
\beq
\| \phi \|_{\hig(S)} \ := \ \sup \{ \| \phi(T) \| \ : \ &T \ {\rm is \ a\ generic\ matrix \ }d{\rm -tuple\ on\ }
\O,
\\
\quad & {\rm and\ } \| \psi(T) \| \leq 1 \ \forall \psi \in S \} .
\eeq

\begin{prop}
\label{propc3}
Let $S$ be a set of holomorphic functions from $\O$ to  ${\D}$.
Then $\hig(S) =  A^\infty(K_S)$ isometrically.
\end{prop}
\bp
Let $\phi$ be in the closed unit ball of $ A^\infty(K_S)$.
Let $T$ be a generic matrix tuple on $\O$, with eigenvectors as in \eqref{eqc3},
and assume that $\| \psi(T) \| \leq 1$ for all $\psi$ in $S$.
Let $F = \{ \lambda_1, \dots, \lambda_n \}$.
Define a kernel $k(z,w)$ on $\O$ by setting it to zero unless both $z$ and $w$ are in $F$,
and on $F$ define 
\[
k(\lambda_i, \lambda_j) \= \la v_i, v_j \ra .
\]
Then $k \in K_S$, so 
\be
\label{eqc4}
(1 - \phi(\lambda_i) \overline{\phi(\lambda_j)})  \la v_i, v_j \ra 
\ \geq \ 0 .
\ee
Then \eqref{eqc4} says that $\| \phi(T) \| \leq 1$, so $\phi$ is in the closed unit ball of
$\hig(S)$.

Conversely, if $\phi$ is in the closed unit ball of
$\hig(S)$, then for every finite set $F \subset \Omega$,
by Theorem~\ref{thmc1} applied to $F$, we have that \eqref{eqc7} holds on $F$.
Hence by the Theorem again, we have $\phi$ is in  the closed unit ball of $ A^\infty(K_S)$.
\ep

\bt
\label{thmc2}
Let $S$ be a set of holomorphic functions from $\O$ to  ${\D}$.
Let $\h$ be an infinite dimensional Hilbert space.
If  $S$ is finite, then
$H^\infty(\h(S)) = A^\infty(K_S)$.
\et
\bp By Theorem~\ref{thm1}, we  have $H^\infty(\h(S))$ is contractively contained
in  $A^\infty(K_S)$. For the converse,
 let $\phi$ be in the closed unit ball of $ A^\infty(K_S)$,
with a network realization formula as above.
Let
$S = \{ \psi_1, \dots, \psi_n \}$.

 Let $\Lambda_j$ be the elements of $C_b(S)$
defined by $\Lambda_j(\psi_i ) = \delta_{ij}$.
Since each $\Lambda_j$ is a projection, we get that  $\rho(\Lambda_j) = P_j$
gives $n$ mutually orthogonal projections that sum to the identity on $\M$.
Then $\rho(E(z)) = \sum_{j=1}^n \psi_j(z) P_j$.

Expanding \eqref{eqc1} as a Neumann series in $ D \rho(E(z))$, the partial sums $\phi_n$ will converge locally uniformly on $\Omega$. Therefore if $T$ is in $\h(S)$, since its spectrum is a compact subset of $\Omega$,
we get that $\phi(T) = \lim_n \phi_n(T)$. We have $\rho(E(T)) =  \sum_{j=1}^n \psi_j(T)\otimes P_j$,
and \eqref{eqc1} extends to
\be
\label{eqc2}
\phi(T) \= {I_\h} \otimes A + ( {I_\h} \otimes B )\rho(E(T)) ( I - ({I_\h} \otimes D) \rho(E(T)))^{-1} {I_\h} \otimes C.
\ee
A calculation with \eqref{eqc2} shows that $ {I_\h} - \phi(T)^* \phi(T) \geq 0$,
so we conclude $\phi \in \SS(\h(S))$.
\ep

\begin{prob}
Let $\h$ be an infinite dimensional Hilbert space.
Do $H^\infty(\h(S))$ and $ A^\infty(K_S)$ coincide
for all non-empty sets $S$ of holomorphic functions
from $\O$ to $\D$?
\end{prob}
              
\bibliographystyle{plain}

\bibliography{../references_uniform_partial}

\end{document}